\documentclass[12pt]{article}

\usepackage{amssymb}
\usepackage{amsmath}
\usepackage{amscd}
\usepackage{ifpdf}

\usepackage[dvips]{graphicx}
\usepackage[ruled]{algorithm}
\usepackage{algorithmic}
\usepackage{palatino, url, multicol}

\numberwithin{equation}{section}

\begin{document}

\title{STORAGE ALLOCATION UNDER PROCESSOR SHARING I: EXACT SOLUTIONS AND ASYMPTOTICS}

\author{
Eunju Sohn\thanks{ Department of Mathematics, Statistics, and
Computer Science, University of Illinois at Chicago, 851 South
Morgan (M/C 249), Chicago, IL 60607-7045, USA. {\em Email:}
esohn3@math.uic.edu.} \and \and Charles Knessl\thanks{ Department of
Mathematics, Statistics, and Computer Science, University of
Illinois at Chicago, 851 South Morgan (M/C 249), Chicago, IL
60607-7045, USA. {\em Email:} knessl@uic.edu.\
\newline\indent\indent{\bf Acknowledgement:} This work was partly supported by NSF grant DMS 05-03745 and NSA grant H 98230-08-1-0102.
}}
\date{ }
\maketitle

\begin{abstract}
\noindent  We consider a processor sharing storage allocation model,
which has $m$ primary holding spaces and infinitely many secondary
ones, and a single processor servicing the stored items (customers). All of the
spaces are numbered and ordered. An arriving customer takes the
lowest available space. We define the traffic intensity $\rho$ to be
$\lambda / \mu$ where $\lambda$ is the customers' arrival rate and $\mu$ is
the service rate of the processor. We study the joint probability
distribution of the numbers of occupied primary and secondary spaces. For $0 < \rho < 1$, we obtain the
exact solutions for $m = 1$ and $m = 2$. For arbitrary $m$ we study the problem in the asymptotic limit $\rho \uparrow 1$ with $m$ fixed. We also develop a semi-numerical semi-analytic method for computing the joint distribution.
\end{abstract}

\setlength{\baselineskip}{24pt}

\section{Introduction}
\indent We consider the following storage allocation model. Suppose that near a restaurant there are $m$ primary parking spaces and across
the street there are infinitely many additional ones. However, the
restaurant has only one waiter who serves all of the customers. All of
the parking spaces are numbered and ordered; the one with rank $= 1$
is closest to the restaurant and the primary spaces are numbered $\{1, 2, 3, ....., m\}$.
  We assume the following: (1) customers arrive
according to a Poisson process with rate $\lambda$, (2) the waiter works at rate $\mu$, (3) an arriving car
parks in the lowest-numbered available space, and (4) if there are $N$
customers in the restaurant, the waiter serves each customer at the
rate $\mu / N$. This corresponds to a processor sharing (PS) service discipline.

\indent Dynamic storage allocation and the fragmentation of computer
memory are among the many applications of this model.  We define
$N_1$ to be the number of occupied primary spaces and $N_2$ to be the number of occupied
secondary spaces.  Then we define $\textbf{S}$ to be the set of the indices of the occupied spaces, and the "wasted spaces" $W$ are defined as the difference between the largest
index of the occupied spaces (Max $\textbf{S}$) and the total
number of occupied spaces ($|\textbf{S}|=N_1+N_2$).  Coffman, Flatto, and Leighton \cite{CFL}
showed that for the processor-sharing model
$$E[W] = \Theta \left( \sqrt{\frac{1}{1-\rho}\log \left(\frac{1}{1-\rho}\right)} \right), \; \; \rho \uparrow 1$$

\noindent where $E[W]$ is the expected value of the wasted spaces.  Here $E[W] = \Theta(f(\rho))$ means that there exist positive constants $c, c'$ such that $c'f(\rho)\leq E[W] \leq c f(\rho)$. Also when $\rho \rightarrow 1$ (the heavy traffic case) Coffman and Mitrani \cite{CM}
obtained upper and lower bounds on $E[W]$ in the form

$$\frac{1}{2}\sqrt{\frac{\pi}{1-\rho}} \leq E[W]\leq\frac{1}{1-\rho}\left(\frac{\pi^2}{6}-1\right).$$

\indent A related model, the $M/M/\infty$ queue with ranked servers, has been studied by many authors \cite{A}, \cite{CK}, \cite{CL}, \cite{K}, \cite{P}, \cite{N}. This differs from the current model in that if there are a total of $N=N_1+N_2$ spaces occupied, the total service rate is $\mu N$, as each customer in the restaurant is served at rate $\mu$. For this model various asymptotic studies appear in \cite{A}, \cite{CL}, \cite{N}. In particular, Aldous \cite{A} showed that the mean number of the wasted spaces is $E[W] \sim \sqrt{2 \rho \log\log \rho}\; $ as $\; \rho = \lambda/\mu \rightarrow \infty$.

\indent A simple derivation of the exact joint distribution of finding $N_1$ (resp., $N_2$) occupied primary (resp., secondary) spaces appears in \cite{S} and detailed asymptotic results for this joint distribution appear in \cite{KB}, \cite{KC}, while the distribution of Max $\textbf{S}$ is analyzed in \cite{KA}. In \cite{KB} Knessl showed how to obtain asymptotic results for the infinite server model directly from the basic difference equation. Since the present processor sharing model does not seem amenable to exact solution, we shall employ such a direct asymptotic approach here.

\indent In this paper, we study the joint probability distribution of the numbers of occupied spaces in the PS model, letting $\pi(k, r) =  Prob[N_1=k, N_2=r]$ in the steady state. From $\pi(k, r)$, we can get the probability distribution of the wasted spaces $W$ from the following relations.
   \begin{eqnarray}
    Prob[W = 0] &=& \sum_{j=0}^\infty \pi(j, 0; j),\label{W0}  \\
    Prob[W = L] &=& \sum_{j=0}^\infty [ \pi(j, 0; L+j) -
    \pi(j, 0; L+j-1)],\;  L \ge 1. \label{WL}
\end{eqnarray}
Then the mean number of the wasted spaces is
\begin{equation*}
    E[W] = \sum_{L=1}^\infty L \; Prob[W=L].
\end{equation*}

\indent The paper is organized as follows. In section 2 we
state the problem and obtain the basic difference equations for $\pi(k,r)$. In section 3 we summarize our main results. Section 4 contains the exact solutions for $m=1$ and $m=2$,  and a sketch of the derivations.  In section 5 we develop a semi-analytic and semi-numerical method for arbitrary $m$, and show that this regains our previous results for $m=1, 2$. Asymptotic solutions for fixed $m$ when $\rho \uparrow 1$ are obtained in section 6. Some numerical studies and comparisons appear in section 7.

\section{Statement of the problem}
We let $N_1 (t)$ (resp., $N_2 (t)$) denote the number of primary
(resp., secondary) spaces occupied at time $t$. The joint steady
state distribution function is
\begin{equation*}
\pi (k, r) = \pi (k, r; m)=\lim_{t\rightarrow\infty} Prob[ N_1 (t) = k,  N_2(t) = r ], \; 0\le k \le m,\; r \ge 0.
\end{equation*}

Let $\rho = \lambda/\mu$ be the traffic intensity and we assume the stability condition $\rho < 1$.
 The pair $(N_1, N_2)$ forms a Markov chain whose transition rates are sketched in \textbf{Fig.1}. The state space is the lattice strip $\{(k,r): 0\leq k \leq m, \;  r\geq 0\}$ and the balance equations are
\begin{multline}\label{bl} (1I_{[k+r>0]}+\rho)\; \pi(k,r)\\ = \rho \; \pi (k-1,
r)I_{[k\geq1]}+\frac{k+1}{k+r+1}\; \pi (k+1, r)I_{[k < m]}\\
 +\frac{r+1}{k+r+1}\; \pi (k, r+1)+\rho\; \pi (m, r-1)I_{[k=m,\; r\geq1]}.
\end{multline}
Here $I$ is an indicator function. The normalization condition is
\begin{equation} \label{n}
\sum_{r=0}^\infty \sum_{k=0}^m \pi(k,r) = 1.
\end{equation}

\indent From our viewpoint we will need to consider explicitly the boundary conditions inherent in (\ref{bl}), so we rewrite the main equation as
\begin{multline} \label{bl2}
(1+\rho)\; \pi(k,r) = \rho \; \pi (k-1,
r)+\frac{k+1}{k+r+1}\; \pi (k+1, r)
 +\frac{r+1}{k+r+1}\; \pi (k, r+1), \\ \; 0\leq k < m, \; r\geq 0, \; k+r > 0.
\end{multline}
and the boundary condition at $k=m$ is
\begin{equation} \label{bcm}
(1+\rho)\; \pi(m,r)\\ = \rho \; \pi (m-1,
r)+\frac{r+1}{m+r+1}\; \pi (m, r+1)+\rho\; \pi (m, r-1), \; r\geq1.
\end{equation}
There are also the two corner conditions
\begin{equation}
\rho \pi(0,0) = \pi(1,0)+\pi(0,1) \label{26}
\end{equation}
and
\begin{equation} \label{27}
(1+\rho)\pi(m,0) = \rho\pi(m-1,0)+\frac{1}{m+1} \pi(m,1).
\end{equation}
In (\ref{bl2}) when $k=0$ we interpret $\pi(-1,r)$ as 0. The boundary condition at $k=m$ in (\ref{bcm}) can be replaced by the
artificial boundary condition
\begin{equation}
\frac{m+1}{m+r+1}\pi(m+1, r) =
\rho\;\pi(m, r-1). \label{ab}
\end{equation}
This is obtained by extending (\ref{bl2}) to hold also at $k=m$ and comparing this to (\ref{bcm}).

\indent We note that the total number $N_1+N_2$ behaves as the number of customers in the $M/M/1-PS$ queue, which is well known to follow a geometric distribution. Thus we have
\begin{equation} \label{29}
\sum_{k+r=N} \pi(k,r) = (1-\rho)\rho^N, \; N \geq 0
\end{equation}
and we can rewrite this as
\begin{eqnarray}
\sum_{r=0}^N \pi(N-r, r) &=& (1-\rho)\rho^N, \; \; 0\leq N \leq m, \\
\sum_{r=N-m}^N \pi(N-r,r) &=& \sum_{k=0}^m \pi(k, N-k) = (1-\rho)\rho^N, \; N \geq m. \label{mm1}
\end{eqnarray}
These identities will provide a useful check on the calculations that follow.
\begin{center}
\includegraphics[angle=0, width=1.0\textwidth]{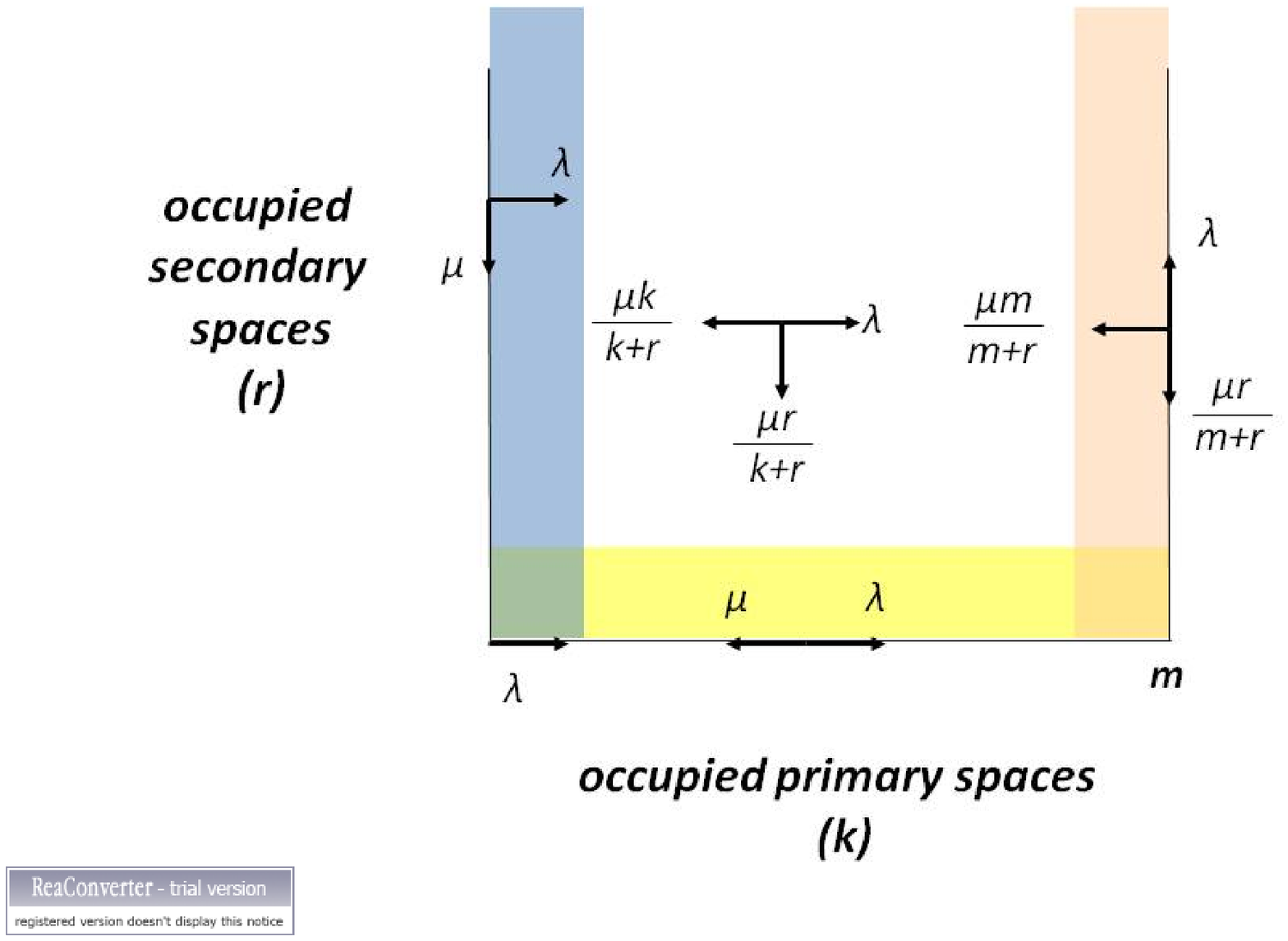}\\
\textbf{Fig. 1}   A sketch of the transition rates.
\end{center}

\section{Summary}
Since the analysis shall become quite involved and technical, we collect here some of the main results.

For $m=1$ the exact solution is given by $\pi(0,0)=1-\rho$ and
\begin{equation}
\pi(0,r) =
(1-\rho)(1+\rho)^{r-1}\int_1^{1+\rho}\left(1-\frac{1}{u}\right)^r
du, \; \; r>0, \label{pi0i3}
\end{equation}
and then $\pi(1,r) = (1-\rho)\rho^{r+1}-\pi(0,r+1).$  For $m=2$ we find that
\begin{eqnarray}
\pi(0,r) &=& (1-\rho)(1+\rho)^{r-2}\rho^{r+2} r(r+1) \int_0^1 \int_0^1 \exp\left[\frac{\rho^2}{(1+\rho)^2} t(1-u)\right] \nonumber \\
&\times& \frac{u^2 (1-u)^{r-1} (1-t)^{r-a^*}}{[1+\rho(1-u)(1-t)]^{r+2}} \; du \; dt,  \label{pi0rm24} \\
\pi(1,r) &=&  (1-\rho)(1+\rho)^{r-1}\rho^{r+2} (r+1)^2 \int_0^1 \int_0^1 \exp\left[\frac{\rho^2}{(1+\rho)^2} t(1-u)\right] \nonumber\\
&\times& u^2 (1-u)^{r-1}  (1-t)^{r-a^*} \frac{[r-2\rho(1-u)(1-t)]}{[1+\rho(1-u)(1-t)]^{r+3}} \; du \; dt.  \label{pi1rm24}
\end{eqnarray}
Here $a^* = \rho/(1+\rho)^2 < 1/4$,  and  $\; \pi(2,r)$ may be computed from $\pi(2,r) = (1-\rho)\rho^{r+2}-\pi(1, r+1)-\pi(0, r+2)$.  Also, $\; \pi(0,0)=1-\rho$,  $\pi(1,0)$ is given above (\ref{pi1rm2}) and then $\pi(0,1)=(1-\rho)\rho-\pi(1,0).$

We next consider the heavy traffic limit where $\varepsilon = 1-\rho \rightarrow 0^+$, with $m=O(1).$  Then setting $r=Y/\varepsilon$, for $Y>0$ and $0\leq k\leq m$ we obtain
\begin{eqnarray} \label{615}
\pi(k,r) = \varepsilon^{m-k+1}\frac{m!}{k!} e^{-Y} Y^{k-m}\left[ 1- \varepsilon \left( \frac{Y}{2} + m + \frac{m^2+(2-k)m-k}{Y}\right)+ O(\varepsilon^2)\right]. \nonumber\\
\end{eqnarray}
This is a two-term asymptotic approximation, but it becomes invalid as $Y\rightarrow 0$.  In section \textbf{6} we discuss the scale $Y=O(\varepsilon)$ (thus $r=O(1)$) where a different approximation to $\pi(k,r)$ must be constructed.

For moderate traffic intensities $\rho<1$ and $m=O(1)$, we have
\begin{equation}
\pi(k,r) \sim (1-\rho)\rho^{m+r}r^{k-m} \frac{m!}{k!}, \; \; r\rightarrow \infty.
\end{equation}
Thus the probability distribution is exponentially small ($O(\rho^{m+r})$), with an additional algebraic decay due to the factor $r^{k-m}$.

\section{Exact solutions for $m=1$ and $m=2$}

We consider $m=1$ and $m=2$ for $0<\rho <1$ fixed.
To solve these problems we use the generating function $F(x, y)$ defined by
\begin{eqnarray}
F(x, y)&=& \sum_{k=0}^{m} \sum_{r=0}^{\infty} \pi (k, r) x^k y^r \label{g1}\\
&=& f_0(y) +xf_1(y)+ \cdot\cdot\cdot\cdot+x^m f_m(y) \label{g2}
\end{eqnarray}
so that
\begin{equation}\label{4}
f_k(y) = \sum_{r=0}^\infty y^r \pi(k,r).
\end{equation}
Multiplying both sides of (\ref{bl}) by $(k+r+1) x^k y^r$ and
summing over $0\leq k\leq m, \; 0\leq r < \infty$ we get a PDE
for $F(x,y)$:
\begin{multline} \label{pdeforF}
(\rho+1-2 \rho x)F + (\rho x -1)(1-x)F_x
+[(\rho+1)y-\rho xy-1]F_y \\ = \rho (y-x)x^m [(m+2)f_m(y)+y f_m'(y)]+1-\rho.
\end{multline}
By setting $x=y$ we find that $F(y,y)=(1-\rho)/(1-\rho y)$.

\indent From (\ref{g1}) and (\ref{g2}) we have $F(x,y) =
f_0(y)+xf_1(y)$ for $m=1$ and $F(x,y)= f_0(y)+xf_1(y)+x^2f_2(y)$ for
$m=2$. Using these forms in (\ref{pdeforF}) and comparing coefficients of
$x^0, x^1$, and (if $m=2$) $x^2, \;$ we obtain the ODEs
$$m=1:$$\begin{eqnarray}
f_1 &=& (1+\rho)f_0+[(1+\rho)y-1]f_0'-1+\rho,   \label{m11} \\
-2\rho f_0-\rho y f_0' &=&[3\rho y-2(1+\rho)]f_1+[\rho y^2-(1+\rho)y+1]f_1'. \label{m12}
\end{eqnarray}
$$m=2:$$\begin{eqnarray}
f_1 &=& (1+\rho)f_0+[(1+\rho)y-1]f_0'-1+\rho,   \label{m21} \\
f_2 &=& (1+\rho)f_1+\frac{(1+\rho)y-1}{2}f_1'-\rho f_0-\frac{\rho}{2} yf_0', \label{m22} \\
-3\rho f_1-\rho yf_1'&=& [4\rho y - 3(1+\rho)] f_2+[\rho y^2-(1+\rho)y+1] f_2'. \label{m23}
\end{eqnarray}
We must thus solve a system of two ODEs if $m=1$ and three ODEs if $m=2$.  We can simplify the calculations by observing that
\begin{equation} \label{F}
F(y,y) = \sum_{k=0}^m \sum_{r=0}^\infty \pi(k,r)y^{k+r} = \sum_{N=0}^\infty y^N \sum_{k+r=N}\pi(k,r) = \frac{1-\rho}{1-\rho y}.
\end{equation}
\noindent Therefore, for $m=1$, $\; f_0(y)+yf_1(y) = (1-\rho)/(1-\rho
y)$ and, for $m=2$, $\; f_0(y)+yf_1(y)+y^2 f_2(y) = (1-\rho)/(1-\rho y)$.  It thus suffices to compute $f_0$ if $m=1$, and $f_0$ and $f_1$ if $m=2$.

\subsection{$m=1$}

By solving the system of equations (\ref{m11}) and (\ref{m12}) for $m=1$, we obtain
\begin{equation} \label{f01}
f_0(y) =
\frac{1-\rho}{[(1+\rho)y-1]^2}\left\{(1+\rho)y^2-(2+\rho)y-y\log
[(1+\rho)(1-\rho y)]+1\right\}
\end{equation}
and
\begin{multline}
f_1(y)=-(1-\rho)+\frac{1-\rho}{(1+\rho)y-1} \\
 \times\left\{(1+\rho)y+\frac{\rho
y}{1-\rho y}+ \frac{\log[(1+\rho)(1-\rho
y)]}{(1+\rho)y-1}-1\right\}.  \label{f11}
\end{multline}
 \noindent Here we also used the
analyticity of $f_0(y)$ at the point $y=1/(1+\rho)$, as $f_0(y)$ must be analytic for $|y| <1$.
We first invert $f_0(y)$ and compute $\pi(0,r)$.  We expand (\ref{f01}) in Taylor series about $y=1/(1+\rho)$, setting $\Delta=y-1/(1+\rho)$, to obtain
\begin{equation} \label{f0s1}
f_0(y) = \frac{1-\rho}{1+\rho} \left\{1+\rho \sum_{n=0}^{\infty}[\rho(1+\rho)\Delta]^n \left(\frac{1}{n+1}+\frac{\rho}{n+2}\right) \right\}.
\end{equation}
Replacing $\Delta$ in (\ref{f0s1}) by $y-1/(1+\rho)$, using a binomial expansion of $\Delta^n$, and reversing the order of summation yields
\begin{equation} \label{f0s2}
f_0(y) = \frac{1-\rho}{1+\rho}\left[1+\rho \sum_{l=0}^{\infty}y^l\sum_{n=l}^{\infty}\rho^n(1+\rho)^l(-1)^{n-l}{n \choose l}\left(\frac{1}{n+1}+\frac{\rho}{n+2}\right)\right].
\end{equation}
From (\ref{f0s2}) we can easily check that $\pi(0,0)=1-\rho$ and then obtain $\pi(0,r)$ for $r>0$ as an infinite series:
\begin{eqnarray}
\pi(0,r) &=& \rho(1-\rho)(1+\rho)^{r-1} \sum_{n=r}^{\infty}\rho^n(-1)^{n-r}{n \choose r}\left(\frac{1}{n+1}+\frac{\rho}{n+2}\right) \nonumber \\
&=& (1-\rho)(1+\rho)^{r-1}\frac{\rho^{r+1}}{r!} \sum_{L=0}^{\infty}(-\rho)^L\frac{(L+r)!}{L!}\left(\frac{1}{L+r+1}+\frac{\rho}{L+r+2}\right). \label{pis2} \nonumber \\
\end{eqnarray}
We can derive another form for $\pi(0,r)$ by expanding (\ref{f01}) around $y=0$ instead of $y=1/(1+\rho)$, which yields
\begin{multline*}
f_0(y) = (1-\rho)\{\; \sum_{h=0}^\infty\sum_{j=0}^\infty(1+\rho)^{h+j}y^{h+j}-[2+\rho+\log(1+\rho)]\sum_{h=0}^\infty\sum_{j=0}^\infty(1+\rho)^{h+j}y^{k+j+1}\\+\sum_{h=0}^\infty\sum_{j=0}^\infty(1+\rho)^{h+j+1}y^{h+j+2}+\sum_{h=0}^\infty\sum_{j=0}^\infty\sum_{l=1}^\infty (1+\rho)^{h+j}\frac{\rho^l}{l}y^{h+j+l+1}\; \}.
\end{multline*}
Therefore, the coefficient of $y^r$ when $r>0$ is
\begin{eqnarray} \label{pis3}
\pi(0,r) &=& (1-\rho)(1+\rho)^{r-1}\{(r+1)(1+\rho)-r[2+\rho+\log(1+\rho)] \nonumber \\
& &+ \;
(r-1)+\sum_{l=1}^{r-1}\frac{r-l}{l}\left(\frac{\rho}{1+\rho}\right)^l\}.
\end{eqnarray}
Furthermore we have
\begin{eqnarray*}
\sum_{l=1}^{r-1}\frac{r-l}{l}\left(\frac{\rho}{1+\rho}\right)^l= r\log(1+\rho)-\rho+\rho^r(1+\rho)^{1-r}-r\sum_{l=r}^{\infty}\frac{1}{l}\left(\frac{\rho}{1+\rho}\right)^l
\end{eqnarray*}
which when used in (\ref{pis3}) gives the following expression for $\pi(0,r)$ when $\; r>0$:
\begin{eqnarray}
\pi(0,r) &=& \left(\frac{1-\rho}{1+\rho}\right)\rho^r\left[1+\rho-\sum_{l=r}^{\infty}\frac{r}{l}\left(\frac{\rho}{1+\rho}\right)^{l-r}\right] \nonumber \\
&=& \left(\frac{1-\rho}{1+\rho}\right)\rho^r\left[\sum_{n=0}^{\infty}\frac{n}{r+n}\left(\frac{\rho}{1+\rho}\right)^{n}\right]. \label{pis4}
\end{eqnarray}
Setting (\ref{pis2}) equal to (\ref{pis4}) leads to the identity
\begin{equation*}
\sum_{n=1}^{\infty}\frac{n}{n+r}\left(\frac{\rho}{1+\rho}\right)^{n} =
 \frac{\rho(1+\rho)^{r}}{r!} \sum_{L=0}^{\infty}(-\rho)^L\frac{(L+r)!}{L!}\left(\frac{1}{L+r+1}+\frac{\rho}{L+r+2}\right),
\end{equation*}
which can be verified directly.  We note that (\ref{pis4}) is much more useful than (\ref{pis2}) for obtaining $\pi(0,r)$ for $r\rightarrow\infty$.

\indent To get $\pi(1,r)$, we can use either (\ref{f11}), or (\ref{pis4}) and the condition $\pi(0,N)+\pi(1, N-1)=(1-\rho)\rho^N, \; N>0$.  We thus obtain
\begin{eqnarray}
\pi(1,r) &=& (1-\rho)\rho^{r+1}-\pi(0,r+1) \label{pi0pi1} \\
&=& \left(\frac{1-\rho}{1+\rho}\right)\rho^{r+1} \left[ \sum_{L=0}^\infty\left(\frac{\rho}{1+\rho}\right)^L-\sum_{L=0}^{\infty}\frac{L}{r+L+1}\left(\frac{\rho}{1+\rho}\right)^{L}\right] \nonumber \\
&=& (1-\rho)\rho^r \sum_{L=1}^\infty \frac{r+1}{r+L}\left(\frac{\rho}{1+\rho}\right)^L. \label{pi1r}
\end{eqnarray}
If $r=0$, from (\ref{pi1r}) we get $\pi(1,0) = (1-\rho)\sum_{L=1}^\infty L^{-1}\rho^{L}(1+\rho)^{-L} = (1-\rho) \log(1+\rho)$.  This is the same as $f_1(0)$ from (\ref{f11}).

\indent We can also derive, from (\ref{pis4}) and (\ref{pi1r}), alternate integral representations for $\pi(0,r)$ and $\pi(1,r)$.  If we write
$1/(r+n)$ in (\ref{pis4}) as $\int_0^1 u^{n+r-1}du$,
(\ref{pis4}) becomes
\begin{eqnarray}
\pi(0,r) &=& \left(\frac{1-\rho}{1+\rho}\right)\rho^r\int_0^1 u^{r-1}\sum_{n=1}^{\infty}n\left(\frac{u \rho}{1+\rho}\right)^{n}du \nonumber \\
&=& (1-\rho)\rho^{r+1}\int_0^1 \frac{u^r}{(1+\rho-u\rho)^2} du, \; r>0. \label{pi0rm1}
\end{eqnarray}
In the same way, we can rewrite (\ref{pi1r}) as
\begin{equation}
\pi(1,r) = (1-\rho)\rho^{r+1}(r+1)\int_0^1 \frac{u^r}{1+\rho-u\rho}du, \; r\geq 0. \label{pi1rm1}
\end{equation}
Other integral representations for $\pi(0,r)$
and $\pi(1,r)$ can be derived from (\ref{f0s1}), by replacing $(n+1)^{-1}$ by $\int_0^1 t^n dt$.
Then $f_0(y)$ becomes
\begin{eqnarray}
f_0(y) &=& \frac{1-\rho}{(1+\rho)^2}\int_0^1 \sum_{n=0}^\infty(1+\rho t)[\rho(\rho+1)\Delta t]^n \rho(\rho+1)dt + \frac{1-\rho}{1+\rho} \nonumber \\
&=& \frac{1-\rho}{1+\rho}\left[\rho \int_0^1 \frac{1+\rho t}{1-\rho(\rho+1)y t+\rho t} dt + 1 \right],  \;  \Delta=y-\frac{1}{1+\rho} \nonumber \\
&=& \frac{1-\rho}{1+\rho}\left[1+ \int_1^{1+\rho} \frac{u}{u-(1+\rho)(u-1)y} du \right],  \;  1+\rho t = u \nonumber \\
&=& \frac{1-\rho}{1+\rho}\left[1+\sum_{n=0}^\infty y^n
(1+\rho)^n\left(\int_1^{1+\rho} \frac{(u-1)^n}{u^n}du\right)\right].
\label{pi0i2}
\end{eqnarray}
From the integral form (\ref{pi0i2}) we can check that
$\pi(0,0)=1-\rho$, and then for $r
> 0$ we obtain (\ref{pi0i3})
and (\ref{pi0pi1}) yields
\begin{eqnarray}
\pi(1,r) &=&
(1-\rho)\rho^{r+1}-(1-\rho)(1+\rho)^{r}\int_1^{1+\rho}\left(1-\frac{1}{u}\right)^{r+1}
du  \nonumber \\
 &=& (1-\rho)(1+\rho)^r \; [\; (1+\rho)\int_1^{1+\rho}\frac{r+1}{u^2}\left(1-\frac{1}{u}\right)^r du \nonumber \\
 & & - \int_1^{1+\rho}\left(1-\frac{1}{u}\right)^{r+1} du \;  ], \label{pi1r3}
\end{eqnarray}
 where we used
 $$\int_{1}^{1+\rho}\frac{r+1}{u^2}\left(1-\frac{1}{u}\right)^r du = \left(\frac{\rho}{1+\rho}\right)^r.$$
 After integrating by parts in the second integral in (\ref{pi1r3}) we obtain
\begin{equation}
\pi(1,r) = (1-\rho)(1+\rho)^r(r+1)\int_1^{1+\rho}\frac{1}{u}\left(1-\frac{1}{u}\right)^r du,  \;  r\geq 0. \label{pi1r2}
\end{equation}

\indent We conclude the analysis of $m=1$ by giving asymptotic
expansions for $\pi(k,r)$ as $r \rightarrow \infty \; $:
\begin{equation*}
\pi(0,r) \sim (1-\rho)\frac{\rho^{r+1}}{r}, \;  \; \;  \pi(1,r) \sim (1-\rho)\rho^{r+1}, \; \; \; r \rightarrow \infty.
\end{equation*}
These follow easily from (\ref{pis4}) and (\ref{pi1r}), or (\ref{pi0rm1}) and (\ref{pi1rm1}), or (\ref{pi0i3}) and (\ref{pi1r2}).

\subsection{$m=2$}
For the $m=2$ case we can use (\ref{m21}) and (\ref{m22}) to express $f_2$ and $f_1$ in terms of $f_0$, and then obtain from (\ref{m23}) the following second order ODE:
\begin{multline}
[(1+\rho)y-1]^2 f_0''(y)+\left[(4\rho^2+7\rho+4)y-\frac{2}{y}-2(1+\rho)\right]f_0'(y)  \\
+ 2
\left[\frac{1}{y^2}+\frac{(1+\rho)}{y}+\rho^2+\rho+1\right]f_0(y) =
2(1-\rho)\left[\frac{1}{y^2(1-\rho y)}+\frac{1}{y}+1+\rho \right].
\label{odem2}
\end{multline}
After some transformations (\ref{odem2}) can be converted into an
inhomogeneous hypergeometric equation.  We can show that both of the
solutions to the homogeneous problem fail to be analytic at
$y=1/(1+\rho)$, and that there is a unique particular solution that
is analytic at this point.  To obtain this solution it is best to
simply expand $f_0(y)$ in Taylor series, by setting
\begin{equation}
f_0 = 1-\rho+y \sum_{n=3}^\infty A_n z^{n-3} \label{f0m2}
\end{equation}
where
\begin{equation}
z=\frac{(1+\rho)^3}{\rho} \Delta =
\frac{(1+\rho)^3}{\rho}\left[y-\frac{1}{1+\rho}\right]. \label{z}
\end{equation}
Using (\ref{f0m2}) in (\ref{odem2}) we can compute the $A_n$
recursively to get
\begin{eqnarray*}
A_n &=& 2(n-1)(n-2)(1-\rho)\frac{(1+\rho)^5}{\rho^3} \sum_{l=n}^\infty
\left(\frac{\rho}{1+\rho}\right)^{2l}\frac{1}{l(l-1)(l-2)}\nonumber
\\& & \times\frac{1}{(n-1-a^*)(n-a^*)\cdot\cdot\cdot(l-2-a^*)(l-1-a^*)}
\end{eqnarray*}
where  $a^* = \rho/(1+\rho)^2. \; $  If we replace $z$ in (\ref{f0m2}) by (\ref{z}) and use the identity $$\Delta^n = \left(y-\frac{1}{1+\rho}\right)^n = \sum_{i=0}^n{n \choose i}\left(\frac{-1}{1+\rho}\right)^{n-i}y^i,$$
we can write $f_0$ as
\begin{eqnarray}
f_0(y) &=& 1-\rho - \frac{2(1-\rho)(1+\rho)^4}{\rho^3} \nonumber\\
&\times&\sum_{i=1}^\infty\sum_{n=i-1}^\infty(n+2)!\;
\Gamma(n+2-a^*)
\left(\frac{1}{a^*}\right)^n\frac{(-1)^{n-i}(1+\rho)^i}{(i-1)!(n-i+1)!}y^i\nonumber\\
&\times&\sum_{l=n+3}^\infty\left(\frac{\rho}{1+\rho}\right)^{2l}\frac{1}{l(l-1)(l-2)}\frac{1}{\Gamma(l-a^*)}.\label{f0m2s}
\end{eqnarray}
Therefore, we find that $\pi(0,0)=1-\rho$ and, for $r>0$,
\begin{eqnarray}
\pi(0,r)&=& - \frac{2(1-\rho)(1+\rho)^4}{\rho^3} \nonumber\\
&\times&\sum_{n=r-1}^\infty(n+2)!\;
\Gamma(n+2-a^*)
\left(\frac{1}{a^*}\right)^{n}\frac{(-1)^{n-r}(1+\rho)^r}{(r-1)!(n-r+1)!} \nonumber \\
&\times&\sum_{l=n+3}^\infty\left(\frac{\rho}{1+\rho}\right)^{2l}\frac{1}{l(l-1)(l-2)}\frac{1}{\Gamma(l-a^*)} \nonumber\\
&=&
\frac{2(1-\rho)(1+\rho)^{3r+2}}{\rho^{r+2}}
\sum_{L=0}^\infty\frac{(L+r+1)!}{L!(r-1)!}\;
\Gamma(L+r+1-a^*) \nonumber \\
&\times&\left(-\frac{1}{a^*}\right)^{L}
\sum_{l=L+r+2}^\infty\left(\frac{\rho}{1+\rho}\right)^{2l}\frac{1}{l(l-1)(l-2)}\frac{1}{\Gamma(l-a^*)} \label{pi0rm2}
\end{eqnarray}
where $L=n-r+1$.

\indent In a similar way we obtain $\pi(1,r)$.  We write $f_1(y)$ as a series in $y$ using (\ref{m21}) and (\ref{f0m2s}) and identify the coefficient of $y^r$, thus obtaining
\begin{eqnarray*}
\pi(1,0) &=& \rho(1-\rho)+\frac{2(1-\rho)(1+\rho)^3}{\rho^2} \; \; [
\; \; \sum_{n=1}^\infty n(n+1)\left(-\frac{1}{a^*}\right)^{n}\Gamma(n+1-a^*)\nonumber\\
&\times&
\sum_{l=n+2}^\infty\left(\frac{\rho}{1+\rho}\right)^{2l}\frac{1}{l(l-1)(l-2)}\frac{1}{\Gamma(l-a^*)}\; \; ]
\end{eqnarray*}
and for $r>0$
\begin{eqnarray}
\pi(1,r) &=&  \frac{2(r+1)}{r!} \frac{(1-\rho)(1+\rho)^{r+5}}{\rho^3} \sum_{n=r-1}^\infty\frac{(n+1)(n+2)!(-1)^{n-r+1}}{(n-r+1)!}
\left(\frac{1}{a^*}\right)^n \nonumber \\
&\times&
\Gamma(n+2-a^*)\sum_{l=n+3}^\infty\left(\frac{\rho}{1+\rho}\right)^{2l}\frac{1}{l(l-1)(l-2)}\frac{1}{\Gamma(l-a^*)}. \label{pi1rm2}
\end{eqnarray}

\indent From (\ref{pi0rm2}) and (\ref{pi1rm2}) we can obtain alternate integral representations for $\pi(0,r)$ and $\pi(1,r)$.  We first rewrite (\ref{pi0rm2}) and (\ref{pi1rm2}) as
\begin{eqnarray}
\pi(0,r) &=& 2(1-\rho)(1+\rho)^{r-2}\rho^{r+2}
\sum_{L=0}^\infty\sum_{j=0}^\infty \frac{(L+r+1)!}{L!(r-1)!}\left(\frac{\rho}{1+\rho}\right)^{2j}\nonumber \\
&\times& (-\rho)^L\frac{\Gamma(L+r+1-a^*)}{\Gamma(j+L+r+2-a^*)} \frac{\Gamma(j+1)}{j!} \nonumber\\
&\times& \frac{1}{(j+L+r+2)(j+L+r+1)(j+L+r)},  \label{pi0rm22}\\
\pi(1,r) &=& 2 (r+1)(1-\rho)(1+\rho)^{r-1}\rho^{r+2} \sum_{L=0}^\infty\sum_{j=0}^\infty\frac{(L+r)(L+r+1)!}{L! \; r!} \nonumber \\
&\times& \left(\frac{\rho}{1+\rho}\right)^{2j}(-\rho)^{L}\frac{\Gamma(L+r+1-a^*)}{\Gamma(j+L+r+2-a^*)}\frac{\Gamma(j+1)}{j!} \nonumber\\
&\times& \frac{1}{(j+L+r+2)(j+L+r+1)(j+L+r)}. \label{pi1rm22}
\end{eqnarray}
We use the beta function
\begin{multline}
B(j+1, L+r+1-a^*) \\
= \frac{\Gamma(j+1)\Gamma(L+r+1-a^*)}{\Gamma(j+L+r+2-a^*)} =
\int_0^1 t^j (1-t)^{L+r-a^*} dt \label{beta}
\end{multline}
and represent the last factors in (\ref{pi0rm22}) and (\ref{pi1rm22}) as the following integral
\begin{eqnarray}
&&\frac{1}{(j+L+r+2)(j+L+r+1)(j+L+r)} \nonumber\\
&&\; \;  = \frac{1}{2(j+L+r+2)}-\frac{1}{j+L+r+1}+\frac{1}{2(j+L+r)}  \nonumber\\
&&\; \;  = \frac{1}{2} \int_0^1 (1-u)^2 u^{j+L+r-1} du  \nonumber \\
&&\; \;  = \frac{1}{2} \int_0^1 u^2 (1-u)^{j+L+r-1} du. \label{frac}
\end{eqnarray}
Using (\ref{beta}) and (\ref{frac}) in (\ref{pi0rm22}) and (\ref{pi1rm22}) leads to
\begin{eqnarray}
\pi(0,r) &=& (1-\rho)(1+\rho)^{r-2}\rho^{r+2} \int_0^1 \int_0^1 u^2 (1-u)^{r-1} (1-t)^{r-a^*}  \nonumber \\
&\times& \sum_{L=0}^\infty\sum_{j=0}^\infty \frac{(L+r+1)!}{L!(r-1)!}(-\rho)^L[(1-u)(1-t)]^L \nonumber\\
&\times& \left(\frac{\rho}{1+\rho}\right)^{2j} \frac{t^j}{j!} (1-u)^j du \; dt, \; \; r>0 \label{pi0rm23}
\end{eqnarray}
and
\begin{eqnarray}
\pi(1,r) &=&  (r+1)(1-\rho)(1+\rho)^{r-1}\rho^{r+2} \int_0^1 \int_0^1 u^2 (1-u)^{r-1}  (1-t)^{r-a^*} \nonumber\\ &\times& \sum_{L=0}^\infty\sum_{j=0}^\infty\frac{(L+r)(L+r+1)!}{L! \; r!} (-\rho)^{L}[(1-u)(1-t)]^L \nonumber \\
&\times& \left(\frac{\rho}{1+\rho}\right)^{2j}\frac{t^j}{j!} (1-u)^j du \; dt, \; \; r>0.  \label{pi1rm23}
\end{eqnarray}
Using the binomial series
\begin{equation}
\sum_{L=0}^\infty \frac{(L+r+1)!}{L!}Z^L = (r+1)!(1-Z)^{-r-2} \label{bin}
\end{equation}
we can explicitly evaluate both infinite series in (\ref{pi0rm23})
and (\ref{pi1rm23}), to obtain the integral
representations (\ref{pi0rm24}) and (\ref{pi1rm24}) for $\pi(0,r)$ and $\pi(1,r)$, for $r>0$.  We can similarly obtain
$\pi(2,r)$ by using (\ref{m22}), but we can also use the identity
$\pi(2,r) = (1-\rho)\rho^{r+2} -\pi(1,r+1)-\pi(0, r+2)$ which holds
for all $r\geq 0$. Also, $\pi(1,0)$ can be computed from
$\rho(1-\rho) - \pi(0,1)$.

\indent We again conclude by giving asymptotic formulas for $\pi(k,r)$ for $r \rightarrow \infty$. These follow easily from the double integrals in (\ref{pi0rm24}) and (\ref{pi1rm24}), since for $r \rightarrow \infty$ the integrands become concentrated in the range $u, t = O(1/r)$. Therefore, evaluating the integrals by the Laplace method we find that as $r\rightarrow \infty$
\begin{eqnarray*}
\pi(0,r) &\sim& \frac{2}{r^2}(1-\rho)\rho^{r+2}, \\
\pi(1,r) &\sim& \frac{2}{r}(1-\rho)\rho^{r+2}, \\
\pi(2,r) &\sim& (1-\rho)\rho^{r+2}.
\end{eqnarray*}

\section{Semi-numerical method for fixed $m$}

In this section we discuss how the exact solutions for fixed $m$ can
be obtained by a semi-numerical approach. We shall reduce the
solution of the two-dimensional problem in (\ref{bl}) to a
one-dimensional one, which must be solved numerically for general
$m$.

\indent If we set $\pi(k,r) = D(N, r)$ where $N=k+r,$  (\ref{bl2}) becomes
\begin{multline}
(N+1)[(1+\rho)D(N,r)-\rho D(N-1,r)-D(N+1,r)] \\= (r+1)D(N+1,r+1)-r D(N+1,r). \label{dbl}
\end{multline}
Since (\ref{dbl}) is separable, we seek solutions in the product form $D(N,r)=\alpha(N)\beta(r)$.  Then (\ref{dbl}) becomes
\begin{equation*}
(N+1)\left[(1+\rho) \; \frac{\alpha(N)}{\alpha(N+1)}-\rho \; \frac{\alpha(N-1)}{\alpha(N+1)}-1\right] = (r+1) \frac{\beta(r+1)}{\beta(r)}-r = -M
\end{equation*}
where $-M$ is a separation constant, which we take to be a negative integer.  Solving the difference equations for $\alpha$ and $\beta$ we obtain
\begin{equation}\label{50}
\alpha(k+r,M) = \frac{1}{2\pi i} \int_C \frac{z^{M-k-r-1}}{(1-z)(1-\rho z)} \; (1-z)^{-M/(1-\rho)}(1-\rho z)^{\rho M/(1-\rho)} dz,
\end{equation}
where  $\; M\leq k+r \;$  and  $\; C \;$ is a small loop around $\; z=0$, and
\begin{equation*}
\beta(r,M) = \frac{(-1)^{M+r}}{r!} \frac{1}{(M-r)!},  \; \;  M \geq r.
\end{equation*}
Therefore, we can write $\pi(k,r)$ as a sum over all $r\leq M\leq k+r$.
\begin{equation}
\pi(k,r) = \sum_{M=r}^{k+r} C(M) \frac{(-1)^{M+r}}{r! (M-r)!} \; \alpha(k+r, M), \; \; r>0 \;  \label{pikrd}
\end{equation}
where $C(M)$ is a function of $M$.

\indent The constants $C(M)$ will be chosen to satisfy the boundary condition along $k=m$, as given by (\ref{bcm}). We are allowing $\pi(k,r)$ to have a different form for $r=0$, and we shall compute $\pi(k,0)$ from (\ref{bl2}) (with $r=0$), (\ref{26}), and (\ref{pikrd}). We define $C(0)$ in (\ref{pikrd}) by requiring (\ref{pikrd}) to hold at $r=0$ if $k=m$.  Thus (\ref{pikrd}) will apply for $\{0\leq k < m, \; r\geq1\}$ or $\{k=m, \; r\geq 0\}$.

\indent When $k=0$, (\ref{50}) yields $\alpha(r,r)=1$ and if $k=1$, we obtain
 $\alpha(r+1, r+1)
= 1$  and  $\alpha(r+1,r) = (1+\rho)(r+1).\; $ Note that $\alpha (N,M) = 0 \; $ if $\; M\geq N$.  Thus (\ref{pikrd}) yields
\begin{equation}
\pi(0,r) = \frac{C(r)}{r!} \label{dpi0r},  \; \; r>0,
\end{equation}
and
\begin{equation}
\pi(1,r) = -\frac{C(r+1)}{r!}+(1+\rho)(1+r)\frac{C(r)}{r!}, \; r>0. \label{dpi1r}
\end{equation}

\noindent If we set $C(r) = r! d(r)$,  (\ref{pikrd}) can be written as
\begin{eqnarray}
\pi(k,r) &=& \sum_{M=r}^{k+r} (-1)^{M-r} d(M) {M \choose r}\alpha(k+r,M)\\
&=& \sum_{l=0}^{k} (-1)^{l}d(r+l){l+r \choose r} A(k, r;l), \label{piA}
\end{eqnarray}
for $r>0$  or  $k=m$,  where $M=r+l$ and
\begin{equation}
A(k,r;l) = \frac{1}{2\pi i} \int_C \frac{z^{l-k-1}}{(1-z)(1-\rho z)} \; (1-z)^{-(l+r)/(1-\rho)}(1-\rho z)^{\rho (l+r)/(1-\rho)} dz. \label{A}
\end{equation}
From (\ref{A})  we obtain
\begin{eqnarray}
A(l,r;l) &=& 1,  \label{al}\\
A(k+1, r-1;l+1) &=& A(k,r;l). \label{al2}
\end{eqnarray}

\indent We solve for $\pi(k,0)$ in terms of $\pi(k,1)$.  If we define $\triangle(k)\equiv \rho^{-k}[\pi(k,0)-\pi(k-1,0)]$, (\ref{bl2}) with $r=0$ can be written as
\begin{equation} \label{5}
\triangle(k)-\triangle(k+1)=\frac{\pi(k,1)}{(k+1)\rho^{k+1}}.
\end{equation}
Summing (\ref{5}) over $1\leq k\leq \hat{k}-1$ gives
\begin{equation}   \label{51}
\triangle(\hat{k}) = \frac{1}{\rho}[\pi(1,0)-(1-\rho)]-\sum_{j=1}^{\hat{k}-1} \frac{\pi(j,1)}{(j+1)\rho^{j+1}}.
\end{equation}
From the corner condition (\ref{26}) and (\ref{51}) we get
\begin{equation*}
\pi(\hat{k},0)-\pi(\hat{k}-1,0)=-\rho^{\hat{k}-1}(1-\rho)^2-\sum_{j=0}^{\hat{k}-1}\rho^{\hat{k}-j-1}\frac{\pi(j,1)}{j+1}.
\end{equation*}
Summing over $1\leq \hat{k}\leq k$ we obtain
\begin{eqnarray}
\pi(k,0) &=& (1-\rho)\rho^k-\sum_{j=0}^{k-1}\sum_{l=0}^{j}\rho^{j-l}\frac{\pi(l,1)}{l+1} \nonumber\\
&=& (1-\rho)\rho^k-\sum_{l=0}^{k-1}\frac{\pi(l,1)}{l+1}\left[\frac{1-\rho^{k-l}}{1-\rho}\right].\label{516}
\end{eqnarray}
This holds also if $k=0$, since then the sum is void and we obtain $\pi(0,0) = 1-\rho$.

We can compute $\pi(m,0)$ (the probability that all $m$ primary and no secondary servers are occupied) from either (\ref{piA}) or (\ref{516}) (with (\ref{piA}) used to compute $\pi(l,1)$ in the sum). By equating these two expressions we find that
\begin{equation}
d(0) = \frac{(1-\rho)\rho^m}{A(m,0;0)}, \label{d0m}
\end{equation}
and also obtain the identity
\begin{equation*}
A(m,0;n) = \sum_{l=n}^{m} \frac{n(1-\rho^{m-l+1})}{l(1-\rho)} A(l,0;n),
\end{equation*}
which can also be derived directly from (\ref{A}).

We still need to determine $d(r)$ (and hence $C(r)$) for $r\geq 1$.  To this end we apply (\ref{bcm}) (or (\ref{ab})), which yields
\begin{eqnarray}
&&\frac{m+1}{m+r+1} \sum_{l=0}^{m+1} (-1)^{l}d(r+l){l+r \choose r} A(m+1, r;l) \nonumber\\
&&= \rho\; \sum_{l=0}^{m} (-1)^{l}d(r+l-1){r+l-1 \choose r-1} A(m, r-1;l), \; r\geq 1 \label{52}
\end{eqnarray}
and then we must also satisfy the corner condition in (\ref{27}).  We can view (\ref{52}) as a difference equation for the $d(r)$, of order $m+2$.  In general this must be solved numerically, though below we discuss the cases $m=1, 2, 3$.  For general $m$ we have thus reduced the two-dimensional problem (\ref{bl2}) to the one-dimensional one in (\ref{52}).

\subsection{$m=1$}
When $m=1$, from (\ref{52}) we obtain the following third order difference equation:
\begin{multline}
2(1+\rho)(r+1)\; d(r+1)-(r+1)\; d(r+2)+\rho(1+\rho)r\; d(r-1) \\ +[\rho(1-r)-(1+\rho)^2(r+1)]\; d(r) = 0, \;  \; r>0. \label{dd}
\end{multline}
Rearranging (\ref{dd}) gives
\begin{multline}\label{53}
\rho(1+\rho)[(r+1)d(r)-r d(r-1)]-(1+2\rho)[(r+1)d(r+1)-r d(r)] \\
+(r+1)d(r+2)-r d(r+1)-[d(r+1)-d(r)]=0, \; \; r>0.
\end{multline}
If we define $\vartheta(r) \equiv \rho(1+\rho)(r+1)d(r)-(1+2\rho)(r+1)d(r+1)+(r+1)d(r+2)-d(r+1)$, then (\ref{53}) yields $\vartheta(r) = \vartheta(r-1)$, so that $\vartheta(r)=\vartheta(0)$ is a constant.
To find this constant we use the corner condition (\ref{27}).  From either (\ref{dpi1r}) (which holds now also at $r=0$) or (\ref{piA}) and (\ref{A}), (\ref{27}) can be rewritten as
\begin{equation}\label{25}
\rho (1-\rho) = (1+\rho)^2 d(0)-2(1+\rho)d(1)+ d(2). \\
\end{equation}
Since $d(0) = \rho (1-\rho)/(1+\rho)$ by (\ref{d0m}), (\ref{25}) implies that
$\vartheta(0)= \rho(1+\rho) d(0)-2(1+\rho)d(1)+d(2)=0$. Thus
$\vartheta(r)= 0$ for all $r\geq 0$ and hence
\begin{equation}
\vartheta(r)=\rho(1+\rho)(r+1)d(r)-[(1+2\rho)(r+1)+1]d(r+1)+(r+1)d(r+2)=0.
\label{const}
\end{equation}
We solve (\ref{const}) by using the generating function
$$G(z) = \sum_{r=0}^\infty  d(r) z^r.$$
Multiplying (\ref{const}) by $z^r$ and summing over all $r$ we obtain
\begin{equation}
z[1-(1+\rho)z]G'-[1+(1+\rho)z]G = -d(0)\frac{1+z}{1-\rho z}.
\label{gode}
\end{equation}
The solution to (\ref{gode}) is, using $d(0) = \rho(1-\rho)/(1+\rho) = G(0)$,
\begin{equation}
G(z) =
\frac{\rho(1-\rho)}{(1+\rho)[1-(1+\rho)z]^2}\left\{1-\frac{1+\rho}{\rho}
\; z \log[(1-\rho z)(1+\rho)]-(1+\rho)z \right\}.   \label{G}
\end{equation}
The integration constant that arises in solving (\ref{gode}) was determined by using the
analyticity of $G(z)$ at $z=1/(1+\rho)$.

We compare $G(z) =
\sum_{r=0}^\infty d(r) z^r $ to $f_0 (z) = \sum_{r=0}^\infty
\pi(0,r)z^r$ in (\ref{4}). Since, for $r>0$, $\pi(0,r) = d(r)$ from (\ref{dpi0r}), the coefficients of $G(z)$ and $f_0(z)$ differ only when $r=0. \; $ Then we use   $$G(0)=d(0)=\rho \frac{\pi(0,0)}{(1+\rho)} = f_0(0) \;  \frac{\rho}{1+\rho}.$$  to find that $f_0(z)-G(z) = (1-\rho)/(1+\rho) = \pi(0,0)-\rho \pi(0,0)/(1+\rho)= \pi(0,0)-d(0)$.

\subsection{$m=2$}
The boundary condition (\ref{bcm}) for $k=m=2$ is
\begin{equation}\label{bckm2}
(1+\rho)\; \pi(2,r) = \rho \; \pi (1,
r)+\frac{r+1}{r+3}\; \pi (2, r+1)+\rho\; \pi (2, r-1).
\end{equation}
From (\ref{piA}) - (\ref{al2})  we obtain
\begin{eqnarray}
\pi(1,r) &=& d(r) A(1,r;0)-d(r+1)(r+1) \nonumber \\
&=& d(r)(1+\rho)(r+1)-d(r+1)(r+1), \label{p1r}
\end{eqnarray}
and
\begin{eqnarray}
\pi(2,r) &=& d(r) A(2, r;0)-d(r+1)(r+1) A(2,r;1)+d(r+2)\frac{(r+1)(r+2)}{2} \nonumber \\
&=& d(r)[(1+\rho+\rho^2)+(1+\rho)^2 r]\frac{(r+2)}{2} \nonumber \\
&-& d(r+1)(1+\rho)(r+1)(r+2)+d(r+2)\frac{(r+1)(r+2)}{2}. \label{p2r}
\end{eqnarray}
Using (\ref{p1r}) and (\ref{p2r}) in (\ref{bckm2}) yields the following fourth order difference equation
\begin{multline}
[d(r+3)(r+1)(r+2)-d(r+2)r(r+1)] \\
-3(1+\rho)[d(r+2)(r+1)(r+2)-d(r+1)r(r+1)]\\
+[d(r+2)r(r+1)-d(r+1)r(r-1)]\\
-\rho(1+\rho+\rho^2)[d(r)(r+2)-d(r-1)(r+1)] \\
-\rho(1+\rho)^2[d(r)r(r+2)-d(r-1)(r-1)(r+1)] \\
+(1+\rho+\rho^2)[d(r+1)(r+3)-d(r)(r+2)]\\
+(1+\rho)^2[d(r+1)(r+1)(r+3)-d(r)r(r+2)]\\
+\rho(1+\rho)[d(r+1)r(r+2)-d(r)(r-1)(r+1)]=0. \label{dm2}
\end{multline}
If we define
\begin{eqnarray}
\vartheta (r) &\equiv&
d(r+3)(r+1)(r+2)-3(1+\rho)d(r+2)(r+1)(r+2)\nonumber\\
&+& d(r+2)r(r+1)-\rho(1+\rho+\rho^2)d(r)(r+2)-\rho(1+\rho)^2d(r)r(r+2)\nonumber\\
&+&(1+\rho+\rho^2)d(r+1)(r+3)+(1+\rho)^2d(r+1)(r+1)(r+3) \nonumber\\
&+&\rho(1+\rho)d(r+1)r(r+2), \label{vartheta}
\end{eqnarray}
then (\ref{dm2}) becomes  $\; \vartheta(r) = \vartheta(r-1)$, so that $\vartheta(r)=\vartheta(0)$ is a constant.
From (\ref{d0m}) we obtain $d(0)$ for $m=2$ as
\begin{equation}
d(0) = \frac{\rho^2 (1-\rho)}{1+\rho +\rho^2}. \label{d02}
\end{equation}
To find the constant $\vartheta(0)$ we again use the corner condition (\ref{27}).
Using (\ref{p2r}), (\ref{d02}), (\ref{516}), and (\ref{piA}), (\ref{27}) for $m=2$ becomes
\begin{equation}
d(3)= \rho^3 (1-\rho)-\frac{1}{2}(6+ 9\rho+6 \rho^2)d(1)
+3(1+\rho)d(2). \label{d3}
\end{equation}
Using (\ref{d3}) in (\ref{vartheta}) we find that
 $\vartheta (0) = 0$ and thus $\vartheta (r)=0 \; $ for
 $r\geq 0$.

Now the coefficient functions in (\ref{vartheta}) are quadratic in $r$, and we introduce again the generating function  $G(z) = \sum_{r=0}^\infty d(r)z^r$ to solve $\vartheta (r)=0$.  This leads to the second order ODE
 \begin{multline}
[1-2z-3\rho z+(1+\rho)(1+3\rho)z^2-\rho(1+\rho)^2 z^3\; ]\; G''  \\
+\left[-\frac{2}{z}-2+(6\rho^2+9\rho+4)z-\rho z^2
(4\rho^2+7\rho+4)\right]G'  \\
+\left[\frac{2}{z^2}+\frac{2}{z}+2-2\rho(\rho^2+\rho+1) z\right]G \\
= 2\frac{\rho^2 (1-\rho)}{1+\rho
+\rho^2}\left(\frac{1}{z^2}+\frac{1}{z}+1\right). \label{odeGm2}
\end{multline}
Since $G(z)= f_0(z)-\pi(0,0)+d(0) =
f_0(z)-(1-\rho)+\rho^2(1-\rho)/(1+\rho+\rho^2)$, if we write
(\ref{odeGm2}) in terms of $f_0$, we obtain (\ref{odem2}) (after dividing (\ref{odeGm2}) by $1-\rho z$).  Equation (\ref{odeGm2}) can be solved as in section 4, to ultimately obtain $d(r)$, and then $\pi(k,r)$ from (\ref{piA}).

\subsection{$m \geq 3$}
\indent This method can be used for general $m$, but we
find that it becomes very difficult to solve the ODE corresponding to (\ref{gode}) or (\ref{odeGm2}) when $m>2$. For
example, when $m=3$ we obtain from (\ref{bcm}) or (\ref{52}) a fifth order difference equation for $d(r)$.  This equation is a perfect difference that may be summed to yield the fourth order difference equation
\begin{multline}\label{thetam3}
(r+1)(r+2)(r+3)d(r+4)+[r(r+1)(r+2)\\
- 4(1+\rho)(r+1)(r+2)(r+3)]d(r+3)+3[(2+\rho)(r+1)(r+2)\\
+ (1+\rho)(1+2\rho)(r+1)(r+2)(r+3)+2(r+1)]d(r+2)\\
+ [3(2+\rho)(r+1)(r+2)-12(1+\rho)(1+\rho+\rho^2)(r+1)(r+2)\\
- (1+\rho)^2(1+4\rho)r(r+1)(r+2)-2(3+2\rho)(2-\rho)(r+1)+6r]d(r+1)\\
+ (1+\rho)[\rho(1+\rho)^2(r-1)r(r+1)+9\rho(1+\rho+\rho^2)r(r+1)\\
- 3\rho(1+\rho^2)(r-2)(r+1)+2\rho^2r]d(r)=0.
\end{multline}
Multiplying (\ref{thetam3}) by $z^r$ and summing over $r$ we
obtain the following third order ODE for $G(z)$ for $m=3$:
\begin{eqnarray}
&&[ \; \rho(1+\rho)^3 z^3-(1+\rho)^2(4\rho+1) z^2+3(1+\rho)(2\rho+1)z-4(1+\rho)\nonumber\\
&&+\left(1+\frac{1}{z}\right) \; ]  \; G''' +3 \; [ \; \rho(1+\rho)(3\rho^2+5\rho+3) z^2 \nonumber\\
&&-(8\rho^3+17\rho^2+13\rho+3)z +2(3 \rho^2+5\rho+3)-\left(1+\frac{1}{z}+\frac{1}{z^2}\right) \; ] \;  G'' \nonumber\\
&&+2 \; [ \; \rho(1+\rho)(9\rho^2+10\rho+9)z -2(1+\rho)(6\rho^2+5\rho+6) \nonumber\\
&&+3\left(1+\frac{1}{z}+\frac{1}{z^2}+\frac{1}{z^3}\right) \; ] \; G'+6 \; [ \; \rho(1+\rho)(1+\rho^2) \nonumber\\
&&-\left(\frac{1}{z}+\frac{1}{z^2}+\frac{1}{z^3}+\frac{1}{z^4}\right) \; ] \; G
=-\frac{6\rho^3(1-\rho)}{(1+\rho)(1+\rho^2)}\left(\frac{1}{z}+\frac{1}{z^2}+\frac{1}{z^3}+\frac{1}{z^4}\right).\nonumber\\
\end{eqnarray}
This appears too difficult to solve in terms of elementary or special functions.  Thus, for general $m$ we shall obtain $\pi(k,r)$ only asymptotically, assuming that $\rho \uparrow 1$.

\section{Asymptotic solutions for $\; \rho \uparrow 1$ with fixed $m$ }
In this section, we fix $m$ and let $\rho \rightarrow 1$.  We thus set $\rho = 1-\varepsilon \; $ where $\varepsilon$ is small and positive.  We shall analyze the problem for two ranges of $r$, $\;  r=Y/\varepsilon = O(\varepsilon^{-1})$ and $r=O(1)$.  Note that necessarily $k\leq m = O(1)$. In this heavy traffic limit most of the probability mass occurs in the range $r=O(\varepsilon^{-1})$.

\subsection{$r = O(\varepsilon^{-1})$}
First we analyze (\ref{bl2}) for $r = O(\varepsilon^{-1})$ by setting
\begin{equation*}
r = \varepsilon^{-1} Y, \; \; \pi(k,r) = \varepsilon^{m-k} {\cal R} (k, Y;\varepsilon).
\end{equation*}
On this scale (\ref{bl2}) becomes
\begin{eqnarray*}
(2-\varepsilon)[Y+\varepsilon(k+1)] {\cal R}(k, Y) &=& \varepsilon(1-\varepsilon)[Y+\varepsilon(k+1)]{\cal R}(k-1,Y) \nonumber\\
+ (k+1){\cal R}(k+1,Y)&+&(Y+\varepsilon){\cal R}(k,Y+\varepsilon),
\end{eqnarray*}
while the boundary condition (\ref{bcm}) at $k=m$ becomes
\begin{eqnarray}\label{623}
(2&-&\varepsilon)[Y+\varepsilon(m+1)] {\cal R}(m, Y) = (Y+\varepsilon){\cal R}(m,Y+\varepsilon) \nonumber\\
&+& (1-\varepsilon)[Y+\varepsilon(m+1)][{\cal R}(m, Y-\varepsilon)+\varepsilon {\cal R}(m-1,Y)].
\end{eqnarray}
We expand ${\cal R}(k,Y)= {\cal R}(k,Y;\varepsilon)$ in the form
\begin{equation*}
{\cal R}(k, Y) = \varepsilon[{\cal R}^{(0)}(k, Y) + \varepsilon
{\cal R}^{(1)}(k, Y) +  O(\varepsilon^2)].
\end{equation*}
From (\ref{bl2}) and (\ref{ab}) we obtain the following equations for ${\cal R}^{(0)}$ and ${\cal R}^{(1)}$, for $Y>0$,
\begin{eqnarray}
Y \; {\cal R}^{(0)}(k,Y) &=& (k+1) \; {\cal R}^{(0)}(k+1,Y),  \;  \;  0\leq k \leq m \label{61}\\
Y \; \frac{\partial}{\partial Y} {\cal R}^{(0)}(k,Y)
&+& (Y-2k-1) \; {\cal R}^{(0)}(k,Y) + Y \; {\cal R}^{(0)}(k-1,Y) \nonumber \\
&=& Y \; {\cal R}^{(1)}(k,Y) - (k+1) \; {\cal R}^{(1)}(k+1,Y),  \; \; 0<k<m. \nonumber\\\label{62}
\end{eqnarray}
The general solution to (\ref{61}) is
\begin{equation}
{\cal R}^{(0)}(k, Y) = \frac{Y^k}{k!} {\cal R}^{(0)}(0,Y), \; \; \;  \; 0\leq k\leq m,  \; Y>0. \label{63}
\end{equation}
Using (\ref{63}) in (\ref{62}) and solving for ${\cal R}^{(1)}$ yields
\begin{equation}
{\cal R}^{(1)}(k, Y) = \frac{Y^k}{k!} {\cal R}^{(1)}(0,Y) -
\frac{Y^{k-1}}{(k-1)!}\left[(Y-1){\cal R}^{(0)}(0,Y)+Y \frac{\partial}{\partial Y} {\cal R}^{(0)}(0,Y)\right]. \label{64}
\end{equation}

It remains to determine ${\cal R}^{(0)}(0,Y)$ and ${\cal R}^{(1)}(0,Y)$.  We use the boundary condition at $k=m$.  From (\ref{623}) we find that this equation is satisfied automatically to orders $O(1)$ and $O(\varepsilon)$, while at orders $O(\varepsilon^2)$ and $O(\varepsilon^3)\; $ we obtain
\begin{equation} \label{65}
\frac{\partial^2}{\partial Y^2} {\cal R}^{(0)}(m,Y) + \frac{\partial}{\partial Y} {\cal R}^{(0)}(m,Y) = 0, \; \; Y>0
\end{equation}
and
\begin{eqnarray}
Y \frac{\partial^2}{\partial Y^2} {\cal R}^{(1)}(m,Y)+(Y-m) \frac{\partial}{\partial Y} {\cal R}^{(1)}(m,Y)
+ Y \frac{\partial}{\partial Y} {\cal R}^{(1)}(m-1,Y)
+ {\cal R}^{(1)}(m-1,Y) \nonumber \\=\left(\frac{Y}{2}-\frac{m}{2}-1\right)\frac{\partial^2}{\partial Y^2} {\cal R}^{(0)}(m,Y)-\frac{Y}{2} \frac{\partial^2}{\partial Y^2} {\cal R}^{(0)}(m-1,Y)
- (m+1)\frac{\partial}{\partial Y} {\cal R}^{(0)}(m,Y)\nonumber \\
-\frac{\partial}{\partial Y} {\cal R}^{(0)}(m-1,Y)-Y \frac{\partial}{\partial Y} {\cal R}^{(0)}(m-2,Y) +{\cal R}^{(0)}(m-1,Y) -{\cal R}^{(0)}(m-2,Y). \nonumber \\\label{66}
\end{eqnarray}
Since ${\cal R}^{(0)}(m,Y) = {\cal R}^{(0)}(0,Y) \; Y^m/m! $  should be integrable about $Y=\infty$, from (\ref{65}) we conclude that ${\cal R}^{(0)}(m,Y)$ is proportional to $e^{-Y}$, so that
\begin{equation}
{\cal R}^{(0)}(0,Y) = \frac{m!}{Y^m} e^{-Y} h(m), \label{67}
\end{equation}
and hence
\begin{equation*}
\pi(k,r) \sim \varepsilon\frac{m!}{k!}\varepsilon^{m-k} Y^{k-m} e^{-Y} h(m)
\end{equation*}
for some function $h(m)$.
The normalization condition (\ref{n}) then gives $h(m) = 1$.
From (\ref{63}) and (\ref{66}) we obtain the following equation for ${\cal R}^{(1)}(0,Y)$
\begin{eqnarray}
Y \frac{\partial^2}{\partial Y^2} {\cal R}^{(1)}(0,Y) &+& (Y+2m) \frac{\partial}{\partial Y} {\cal R}^{(1)}(0,Y) + \left(\frac{m^2}{Y}-\frac{m}{Y}+m\right) {\cal R}^{(1)}(0,Y)    \nonumber \\
&=& -\left(\frac{2m^2+4m}{Y^2}+\frac{m^2+2m}{Y}-\frac{Y}{2}\right){\cal R}^{(0)}(0,Y) \nonumber\\
&=&  -\left(\frac{2m^2+4m}{Y^2}+\frac{m^2+2m}{Y}-\frac{Y}{2}\right) \frac{m!}{Y^m} e^{-Y}.  \label{69}
\end{eqnarray}
Solving the ODE (\ref{69}) we find that the general solution takes the form
\begin{equation} \label{610}
{\cal R}^{(1)}(0,Y) = C_1(m) Y^{-m}+ C_2(m) Y^{-m} e^{-Y} -\left(m^2+2m +\frac{Y^2}{2}\right)\frac{m!}{Y^{m+1}} e^{-Y}.
\end{equation}
To fully determine ${\cal R}^{(1)}(0,Y)$ we use (\ref{mm1}), which can be written as, for $N\geq m$,
\begin{eqnarray}
\varepsilon(1-\varepsilon)^N &=& \sum_{L=0}^m \pi(m-L, N+L-m) \nonumber \\
 &=& \sum_{L=0}^m \varepsilon^L {\cal R}(m-L, \varepsilon N+\varepsilon L-\varepsilon m). \label{611}
\end{eqnarray}
Setting $N=Z/\varepsilon$ and using
\begin{equation*}
(1-\varepsilon)^{Z/\varepsilon} = \exp\left[\frac{Z}{\varepsilon}\log(1-\varepsilon)\right] = e^{-Z}\left[1-\frac{\varepsilon}{2}Z+O(\varepsilon^2)\right],
\end{equation*}
from (\ref{611}) we find that
\begin{eqnarray} \label{613}
{\cal R}^{(0)}(m,Z)&+&\varepsilon\left[{\cal R}^{(1)}(m,Z)-m \frac{\partial}{\partial Z} {\cal R}^{(0)}(m,Z) +{\cal R}^{(0)}(m-1,Z)\right] \nonumber \\
&=& e^{-Z}[1-\frac{\varepsilon}{2}Z+O(\varepsilon^2)].
\end{eqnarray}
By comparing terms for order $O(1)$ in (\ref{613}) we obtain ${\cal R}^{(0)}(m,Z) = e^{-Z}$, which is the same as what we obtained from (\ref{63}) and (\ref{67}).  Using (\ref{63}), (\ref{64}), (\ref{67}) and (\ref{613}), the $O(\varepsilon)$ terms in (\ref{613}) yield
\begin{equation}\label{614}
{\cal R}^{(1)}(0,Z) = -\frac{m!}{Z^m} e^{-Z} \varepsilon\left[\frac{Z}{2}+m+ \frac{m^2+2m}{Z}\right].
\end{equation}
Therefore, by comparing (\ref{610}) and (\ref{614}) we find that $C_1=0$ and $C_2= - m \;  m!$ in (\ref{610}), and thus obtain (\ref{615}),
which is a two-term asymptotic approximation to $\pi(k,r)$ for $r=O(\varepsilon^{-1})$.  Note that the second term becomes comparable to the leading term for small $Y$, with $Y=O(\varepsilon)$.

\subsection{$r=O(1)$}
We next examine the problem for $r=O(1)$. Now both $k$ and $r$ are $O(1)$ and we set  $\pi(k,r)=\varepsilon Q(k,r;\varepsilon)$, which satisfies
\begin{eqnarray} \label{616}
(2-\varepsilon)(k+r+1) Q(k, r) &=& (k+1)Q(k+1,r)+(r+1)Q(k,r+1)\nonumber \\
&+& (1-\varepsilon)(k+r+1)Q(k-1,r)
\end{eqnarray} for $0<k<m, \; r\geq 0$,
with the boundary condition
\begin{eqnarray}\label{619}
(2-\varepsilon) Q(m,r) &=& (1-\varepsilon)  Q (m-1,r)+(1-\varepsilon) Q (m, r-1) \nonumber\\
&+&\frac{r+1}{m+r+1}\; Q (m, r+1), \; r\geq1,
\end{eqnarray}
and the corner conditions
\begin{equation*}
(2-\varepsilon)(m+1) Q(m,0) = (1-\varepsilon)(m+1)  Q (m-1,0)+ Q (m, 1),
\end{equation*}
and
\begin{equation}\label{624}
(1-\varepsilon)Q(0,0)=Q(1,0)+Q(0,1).
\end{equation}

We expand $Q(k,r)$ in the form
$$Q(k,r) =  Q^{(0)}(k,r) + \varepsilon Q^{(1)}(k,r)+O(\varepsilon^2).$$  Then from (\ref{616}) we obtain at the first two orders
\begin{equation}  \label{617}
2\; Q^{(0)}(k,r) - Q^{(0)} (k-1,r) = \frac{k+1}{k+r+1}\; Q^{(0)} (k+1, r)
 +\frac{r+1}{k+r+1}\; Q^{(0)} (k, r+1),
\end{equation}
and
\begin{eqnarray}
 2\; Q^{(1)}(k,r)- Q^{(1)} (k-1,r)&-&\frac{k+1}{k+r+1}\; Q^{(1)} (k+1, r)-\frac{r+1}{k+r+1}\; Q^{(1)} (k, r+1)\nonumber\\
&=& Q^{(0)}(k,r)-Q^{(0)}(k-1,r). \label{618}
\end{eqnarray}
Equations (\ref{617}) and (\ref{618}) hold for $\; 0\leq k < m, \; r\geq 0, \; k+r > 0,$  with $Q(-1,r)=0$.
We also obtain boundary conditions at $k=m$ from (\ref{619}) as
\begin{equation}\label{621}
2Q^{(0)}(m,r)-Q^{(0)}(m-1,r)-Q^{(0)}(m,r-1) = \frac{r+1}{m+r+1} Q^{(0)}(m,r+1)
\end{equation}
and
\begin{eqnarray*}
 2Q^{(1)}(m,r)&-&Q^{(1)}(m-1,r)-Q^{(1)}(m,r-1)-\frac{r+1}{m+r+1} Q^{(1)}(m,r+1)  \nonumber \\
&=& Q^{(0)}(m,r)-Q^{(0)}(m-1,r)-Q^{(0)}(m,r-1).
\end{eqnarray*}
The corner conditions at $(k,r)=(0,0)$ are given by, in view of (\ref{624}),
\begin{eqnarray}
Q^{(0)}(1,0)+Q^{(0)}(0,1) &=& 1  \label{622}\\
Q^{(1)}(1,0)+Q^{(1)}(0,1) &=& -1
\end{eqnarray}
where we used $Q^{(0)}(0,0)=1$ and $Q^{(1)}(0,0)=0$.

We also require that the expansions on the $Y$ and $r$ scales asymptotically match, in an intermediate limit where $Y \rightarrow 0$ and $r \rightarrow \infty$.  This means that
\begin{multline} \label{620}
\varepsilon^{m-k}\; [\; {\cal R}^{(0)}(k, Y) + \varepsilon
{\cal R}^{(1)}(k, Y) + O(\varepsilon^2) \; ] \mid_{Y\rightarrow 0} \\
\sim \;  Q^{(0)}(k,r) + \varepsilon Q^{(1)}(k,r)+O(\varepsilon^2) \mid_{r \rightarrow \infty}.
\end{multline}
By setting $Y=\varepsilon r$ in the left side of (\ref{620}) and using (\ref{615}) we find that, for $r \rightarrow \infty$,
\begin{equation}\label{627}
Q^{(0)}(k,r) = \frac{m!}{k!} \; r^{k-m} \left[1- \frac{m^2+(2-k)m-k}{r}+O(r^{-2})\right]
\end{equation}
and
\begin{equation}\label{628}
Q^{(1)}(k,r) = \frac{m!}{k!} \; r^{k-m}[-r+(m-k)(m+1)+O(r^{-1})].
\end{equation}
We also get $\sum_{k+r=N} Q^{(0)}(k,r) = 1$  and $\sum_{k+r=N} Q^{(1)}(k,r) = -N$  from (\ref{29}).

The basic problem (\ref{617}), (\ref{621}), and (\ref{622}) for the leading term $Q^{(0)}(k,r)$ is only slightly simpler than the full problem in (\ref{616}) (or (\ref{bl2})).  The only simplification is that $\rho$ has been replaced by 1.  We have not been able to solve this for general $m$, but the results in section 4 can be used to identify $Q^{(0)}(k,r)$ and $Q^{(1)}(k,r)$ when $m=1$ and $m=2$.  For $m=1$ we use (\ref{pi0i3})  and (\ref{pi1r3}) and expand these expressions for $\rho \rightarrow 1$ and $r>0$, to find that
\begin{eqnarray}
Q^{(0)}(0,r) &=& 2^{r-1}\int_1^2 \left(1-\frac{1}{u}\right)^r \; du, \label{625}\\
Q^{(1)}(0,r) &=& -(r-1) 2^{r-2}\int_1^2\left(1-\frac{1}{u}\right)^r \; du -\frac{1}{2}, \\
Q^{(0)}(1,r) &=& 1-2^r\int_1^2\left(1-\frac{1}{u}\right)^{r+1} \; du \\
&=& 1-Q^{(0)}(0,r+1), \nonumber\\
Q^{(1)}(1,r) &=& -r-\frac{1}{2}+2^{r-1}r\int_1^2 \left(1-\frac{1}{u}\right)^{r+1} \; du \\
&=& -r-1-Q^{(1)}(0,r+1). \nonumber
\end{eqnarray}
For $m=2 \; $ (\ref{pi0rm24}) and (\ref{pi1rm24}) give the leading terms as, for $r>0$,
\begin{eqnarray}
Q^{(0)}(0,r) &=& 2^{r-2} r(r+1) \int_0^1 \int_0^1 \exp\left[\frac{t(1-u)}{4}\right] \frac{u^2 (1-u)^{r-1} (1-t)^{r-1/4}}{(2-t-u+ut)^{r+2}} \; du \; dt \nonumber \\
 \\
Q^{(0)}(1,r) &=& 2^{r-1} (r+1)^2 \int_0^1 \int_0^1 \exp\left[\frac{t(1-u)}{4}\right] u^2 (1-u)^{r-1} (1-t)^{r-1/4}\nonumber\\
&\times&  \frac{r-2(1-u)(1-t)}{(2-t-u+ut)^{r+3}} \; du \; dt.\label{629}
\end{eqnarray}
and $Q^{(0)}(2,r)$ can be computed from
\begin{equation}\label{626}
Q^{(0)}(2,r) = 1-Q^{(0)}(1, r+1)-Q^{(0)}(0,r+2).
\end{equation}

The results of section 4 can also be used to identify the correction term $Q^{(1)}(k,r)$ for $m=2$ and $k=0, 1, 2$.  From (\ref{625})-(\ref{626}) we can easily verify that (\ref{627}) and (\ref{628}) are satisfied for $m=1$ and $m=2$, by expanding the integrals in (\ref{625})-(\ref{629}) for $r\rightarrow \infty$ by the Laplace method.

\section{Numerical Studies}

We assess the accuracy of some of the asymptotic formulas we obtained.  First we consider the heavy traffic case, where (\ref{615}) applies.  In \textbf{Table 1} we consider $\varepsilon = 1-\rho = 0.1, 0.05,  0.02$ and $0.01$.  We take $m=3$ and $Y=\varepsilon r =1$, and compare the one and two term approximations in (\ref{615}) to the exact (numerical) values of $\pi(k, r)$.  For each value of $\varepsilon$ we give $\pi(k, r)$ for $0\leq k\leq 3$.  The two-term approximation includes the $O(\varepsilon)$ term inside the brackets in (\ref{615}).
When $\varepsilon =0.1$ we see that the one-term approximation is quite poor, while two-term approximation is even worse, and may leads to a negative answer.  However, as $\varepsilon$ decreases to 0.01 the agreement becomes quite good, and we also clearly see the improvement obtained by using the second term in (\ref{615}).

In \textbf{Table 2} we consider the limit $r\rightarrow \infty$ with $k, m = O(1)$, where we have the asymptotic formula $\pi(k, r) \sim (1-\rho)\rho^{m+r}r^{k-m} m!/k!$.  We begin with $r=5$, which we ultimately increase to $r=50$.  The agreement between exact and asymptotic results is quite poor when $r=5$ but becomes much better for $r=50$.  For each value of $r$ the agreement is the worst when $k=0$ and improves with increasing $k$.  These data suggest that it may be useful to compute also the correction term, as we did in (\ref{615}).

These comparisons show that the asymptotics agree  reasonably well with the exact numerical values of $\pi(k, r)$.  In some cases it proves useful to compute more than one term in the asymptotic series.

\newpage

\begin{center}
\textbf{Table 1}
\begin{equation*}
\rho =1-\varepsilon, \;  \;  m=3, \;  \;  Y=1,  \;  \;  r=Y/\varepsilon.
\end{equation*}
\end{center}

\begin{center}
 \begin{tabular} {|c|c||c|c|c|}       \hline
$\varepsilon$ &$k$     & exact & one-term & two-term  \\
\hline $0.1$ & $0$ & $5.40\times10^{-5}$ & $.000220$ & $<0$\\
\hline $   $ & $1$ & $.000696$           & $.00220$  & $<0$\\
\hline $   $ & $2$ & $.00463$            & $.0110$   & $<0$\\
\hline $   $ & $3$ & $.0210$             & $.0367$   & $.0128$\\
\hline $0.05$ & $0$ & $6.28\times10^{-6}$ & $1.37\times10^{-5}$ & $1.03\times10^{-6}$\\
\hline $   $  & $1$ & $.000146$      & $.000275$                & $7.58\times10^{-5}$\\
\hline $   $  & $2$ & $.00172$       & $.00275$                 & $.00131$\\
\hline $   $  & $3$ & $.0136$        & $.0183$                  & $.0124$\\
\hline $0.02$ & $0$ & $2.50\times10^{-7}$ & $3.53\times10^{-7}$ & $2.22\times10^{-7}$\\
\hline $   $  & $1$ & $1.34\times10^{-5}$ & $1.76\times10^{-5}$ & $1.25\times10^{-5}$\\
\hline $   $  & $2$ & $.000361$           & $.000441$           & $.000348$\\
\hline $   $  & $3$ & $.00649$            & $.00735$            & $.00640$\\
\hline $0.01$ & $0$ & $1.84\times10^{-8}$ & $2.20\times10^{-8}$ & $1.79\times10^{-8}$\\
\hline $   $  & $1$ & $1.91\times10^{-6}$ & $2.20\times10^{-6}$ & $1.88\times10^{-6}$\\
\hline $   $  & $2$ & $9.96\times10^{-5}$ & $.000110$           & $9.87\times10^{-5}$\\
\hline $   $  & $3$ & $.00345$            & $.00367$            & $.00343$\\
\hline
\end{tabular}\\
\end{center}

\newpage

\begin{center}
\textbf{Table 2}
\begin{equation*}
\rho =0.5, \;  \;  m=3,  \; \; 5\leq r\leq 50.
\end{equation*}
\end{center}

\begin{center}
 \begin{tabular} {|c|c||c|c|c|}       \hline
 $r$        & $k$     & exact & asymptotic  \\
\hline $ 5$ & $0$ & $2.29\times10^{-5}$ & $9.38\times10^{-5}$\\
\hline $  $ & $1$ & $1.58\times10^{-4}$ & $4.69\times10^{-4}$\\
\hline $  $ & $2$ & $6.02\times10^{-4}$ & $1.17\times10^{-3}$\\
\hline $  $ & $3$ & $1.65\times10^{-3}$ & $1.95\times10^{-3}$\\
\hline $10$ & $0$ & $1.60\times10^{-7}$ & $3.66\times10^{-7}$\\
\hline $  $ & $1$ & $1.94\times10^{-6}$ & $3.66\times10^{-6}$\\
\hline $  $ & $2$ & $1.23\times10^{-5}$ & $1.83\times10^{-5}$\\
\hline $  $ & $3$ & $5.49\times10^{-5}$ & $6.10\times10^{-5}$\\
\hline $20$ & $0$ & $2.83\times10^{-11}$ & $4.47\times10^{-11}$\\
\hline $  $ & $1$ & $6.29\times10^{-10}$ & $8.94\times10^{-10}$\\
\hline $  $ & $2$ & $7.18\times10^{-9}$  & $8.94\times10^{-9}$\\
\hline $  $ & $3$ & $5.60\times10^{-8}$  & $5.96\times10^{-8}$\\
\hline $30$ & $0$ & $9.41\times10^{-15}$ & $1.29\times10^{-14}$\\
\hline $  $ & $1$ & $3.04\times10^{-13}$ & $3.88\times10^{-13}$\\
\hline $  $ & $2$ & $5.00\times10^{-12}$ & $5.82\times10^{-12}$\\
\hline $  $ & $3$ & $5.57\times10^{-11}$ & $5.82\times10^{-11}$\\
\hline $40$ & $0$ & $4.18\times10^{-18}$ & $5.33\times10^{-18}$\\
\hline $  $ & $1$ & $1.77\times10^{-16}$ & $2.13\times10^{-16}$\\
\hline $  $ & $2$ & $3.79\times10^{-15}$ & $4.26\times10^{-15}$\\
\hline $  $ & $3$ & $5.49\times10^{-14}$ & $5.68\times10^{-14}$\\
\hline $50$ & $0$ & $2.19\times10^{-21}$ & $2.66\times10^{-21}$\\
\hline $  $ & $1$ & $1.15\times10^{-19}$ & $1.33\times10^{-19}$\\
\hline $  $ & $2$ & $3.03\times10^{-18}$ & $3.31\times10^{-18}$\\
\hline $  $ & $3$ & $5.40\times10^{-17}$ & $5.55\times10^{-17}$\\
\hline
\end{tabular}\\
\end{center}

\end{document}